\documentclass[12pt,a4paper]{article}

\usepackage[utf8]{inputenc} 
\usepackage{geometry} 
\usepackage{graphicx} 

\usepackage[colorlinks,citecolor=red]{hyperref}
\usepackage{amsmath,amsthm}
\usepackage{amsfonts}
\usepackage{dsfont}
\usepackage{mathrsfs}
\usepackage{latexsym}
\usepackage{graphicx}
\usepackage{amssymb}
\usepackage{float}
\usepackage{epsfig}
\usepackage{epstopdf}
\usepackage{color,xcolor}
\usepackage{booktabs} 
\usepackage{array} 
\usepackage{paralist} 
\usepackage{verbatim} 
\usepackage{subfig} 

\newtheorem{theorem}{Theorem}[section]
\newtheorem{lemma}[theorem]{Lemma}

\newtheorem{problem}[theorem]{Problem}

\newtheorem{conjecture}[theorem]{Conjecture}

\usepackage{fancyhdr} 
\usepackage{sectsty}
\allsectionsfont{\sffamily\mdseries\upshape} 

\usepackage[nottoc,notlof,notlot]{tocbibind} 
\usepackage[titles,subfigure]{tocloft} 

\title{{\bf A note on the $A_{\alpha}$-spectral radius of graphs}\thanks{Supported by National Natural Science Foundation of China (Nos. 11401211 and 11471121)
and Fundamental Research Funds for the Central
Universities (No. 222201714049).}}
\author{ Huiqiu Lin$^{a}$\thanks{Corresponding author. Email:~huiqiulin@126.com (H. Lin).},~~Xing Huang$^{a}$,~~Jie Xue$^{b}$
\\
{\footnotesize $^a$Department of Mathematics, East China University of Science and Technology, Shanghai, PR China}\\
{\footnotesize $^b$Department of Computer Science and Technology, East China Normal University, Shanghai, PR China }}

\date{}

\begin{document}
\maketitle

\begin{abstract}
Let $G$ be a graph with adjacency matrix $A(G)$ and let $D(G)$
be the diagonal matrix of the degrees of $G$. For any real $\alpha\in [0,1]$, Nikiforov [Merging the $A$- and $Q$-spectral theories, Appl. Anal. Discrete Math. 11 (2017) 81--107] defined the matrix $A_{\alpha}(G)$ as
$A_{\alpha}(G)=\alpha D(G)+(1-\alpha)A(G).$
Let $u$ and $v$ be two vertices of a connected graph $G$. Suppose that $u$ and $v$ are connected by a path $w_0(=v)w_1\cdots w_{s-1}w_s(=u)$ where $d(w_i)=2$ for $1\leq i\leq s-1$.
Let $G_{p,s,q}(u,v)$ be the graph obtained by attaching the paths $P_p$ to $u$ and $P_q$ to $v$.
Let $s=0,1$. Nikiforov and Rojo [On the $\alpha$-index of graphs with pendent paths, Linear Algebra Appl. 550 (2018) 87--104] conjectured that $\rho_{\alpha}(G_{p,s,q}(u,v))<\rho_{\alpha}(G_{p-1,s,q+1}(u,v))$ if $p\geq q+2.$
In this paper, we confirm the conjecture. As applications, firstly, the extremal graph with maximal $A_{\alpha}$-spectral radius with fixed order and cut vertices
is characterized. Secondly, we characterize the extremal tree which attains the maximal $A_{\alpha}$-spectral radius
with fixed order and matching number. These results generalize some known results.


\end{abstract}

\section{Introduction}
All graphs considered here are simple and undirected. Let $G$ be a graph with adjacency matrix $A(G)$, and let $D(G)$ be the diagonal matrix of the degrees of $G$. For any real $\alpha\in [0,1]$, Nikiforov \cite{VN1} defined the matrix $A_{\alpha}(G)$ as
$$A_{\alpha}(G)=\alpha D(G)+(1-\alpha)A(G).$$
It is clear that $A_{\alpha}(G)$ is the adjacency matrix if $\alpha=0$, and $A_{\alpha}(G)$ is essentially equivalent to signless Laplacain matrix if $\alpha=1/2$.
Write $\rho_{\alpha}(G)$ for the spectral radius of $A_{\alpha}(G)$ and call it $A_{\alpha}$-spectral radius of $G$ (or the $\alpha$-\emph{index }of $G$).

Let $G$ be a connected graph and $u,v$ be two distinct vertices of $V(G)$.
Let $G_{p,q}(u,v)$ be the graph obtained by attaching the paths $P_p$ to $u$ and $P_q$ to $v$.
Nikiforov and Rojo \cite{NO} posed the following problem which is inspired by the result of Li and Feng \cite{LF}.
\begin{problem}\rm{(\cite{NO}, Problem 20)}\label{pro1}
For which connected graphs $G$ the following statement is true:
Let $\alpha\in[0,1)$ and let $u$ and $v$ be non-adjacent vertices of $G$ of degree at least 2.
If $q\geq1$ and $p \geq q+2$, then
$\rho_{\alpha}(G_{p,q}(u,v))<\rho_{\alpha}(G_{p-1,q+1}(u,v))$.
\end{problem}
Let $G$ be a connected graph and $u,v\in V(G)$ with $d(u),d(v)\geq 2$. Suppose that $u$ and $v$ is connected by a path $w_0(=v)w_1\cdots w_{s-1}w_s(=u)$ where $d(w_i)=2$ for $1\leq i\leq s-1$.
Let $G_{p,s,q}(u,v)$ be the graph obtained by attaching the paths $P_p$ to $u$ and $P_q$ to $v$.
In particular, $s=0$ means that $u$ and $v$ are the same vertex and so the graph $G_{p,0,q}(v,v)$ is the graph obtained by attaching the paths $P_p$ and $P_q$ to $v$.
We first show Problem \ref{pro1} is true for the type of graphs $G_{p,s,q}(u,v)$ when $p-q\geq \max\{s+1,2\}$.
\begin{theorem}\label{thm1}
Let $0\leq \alpha <1$. If $p-q\geq \max\{s+1,2\}$, then $\rho_{\alpha}(G_{p-1,s,q+1}(u,v))>\rho_{\alpha}(G_{p,s,q}(u,v))$.
\end{theorem}

When $s=0,1$, Theorem \ref{thm1} implies the following conjecture is true which is also posed in \cite{NO}.
\begin{conjecture}\rm{(\cite{NO}, Conjectures 18 and 19)}\label{conj1}
Let $0\leq \alpha<1$ and $s=0, 1$. If $p\geq q+2$, then $\rho_{\alpha}(G_{p,s,q}(u,v))<\rho_{\alpha}(G_{p-1,s,q+1}(u,v))$.
\end{conjecture}
It needs to notice that Conjecture \ref{conj1} was independently confirmed by Guo and Zhou \cite{GZ}.
Brualdi and Solheid \cite{BS} proposed the following problem concerning the spectral radius of
graphs:
\begin{problem}
Given a set of graphs, to find an upper bound for the spectral radius in this set and
characterize the graphs in which the maximal spectral radius is attained.
\end{problem}

Nikiforov \cite{VN1} determined the graph with maximal $A_{\alpha}$-spectral radius among graphs with fixed order and chromatic number.
The graph with maximal $A_{\alpha}$-spectral radius among all $K_{r+1}$-free graphs was also determined.
In \cite{VN2}, it was presented that the star and the path attain the maximal and the minimal $A_{\alpha}$-spectral radius among all trees, respectively.
Note also that the path attains the minimal $A_{\alpha}$-spectral radius for all connected graphs.
Very recently, the graph with maximal $A_{\alpha}$-spectral radius among all graphs with given order and diameter
and the graph with minimal $A_{\alpha}$-spectral radius among all graphs with given order and clique number were considered in
\cite{NO,X}.
In this paper, we will characterize some extremal graphs with maximal $A_{\alpha}$-spectral radius among some other type of graphs.

A \emph{cut vertex} in a connected graph $G$ is a vertex whose deletion breaks
the graph into two (or more) parts. Let $G_{n,k}$ be the graph obtained
from the complete graph $K_{n-k}$ by attaching paths of almost equal
lengths to all vertices of $K_{n-k}$. Berman and Zhang \cite{BZ} and Zhu \cite{BZhu} showed that $G_{n,k}$ attains the maximal spectral radius, signless Laplacian spectral radius
among all graphs with order $n$ and $k$ cut vertices, respectively. We generalize their result to $A_{\alpha}$-spectral radius as follows.
\begin{theorem}\label{thm2}
Let $G$ be a connected graph with order $n$ and $k$ cut vertices. If $0\leq\alpha<1$, then $\rho_{\alpha}(G)\leq \rho_{\alpha}(G_{n,k})$,
with equality if and only if $G\cong G_{n,k}.$
\end{theorem}

Two distinct edges in a graph $G$ are \emph{independent} if they are not incident with a common vertex
in $G$. A set of pairwise independent edges in $G$ is called a \emph{matching} in $G$. A matching of maximum
cardinality is a \emph{maximum matching} in $G$. The \emph{matching number} $m(G)$ of $G$
is the cardinality of a maximum matching of $G$. It is well known that $m(G)\leq n/2$
with equality if and only if $G$ has a perfect matching. A vertex $v$ is \emph{matched} if it is incident to an edge in the matching; otherwise the vertex is \emph{unmatched}.
We define a tree $A(n,k)$, $n\geq 2k$, with $n$ vertices as follows: $A(n,k)$ is obtained from the star
graph $S_{n-k+1}$ by attaching a pendent edge to each of certain $k-1$ non-central vertices of $S_{n-k+1}$. Note that $A(n,k)$ has an $k$-matching.
Hou and Li \cite{HL} and Guo \cite{Guo} proved that $A(n,k)$, $n\geq 2k$, attains the maximal spectral radius, signless Laplacian spectral radius
among all graphs with order $n$ and matching number $k$, respectively. We generalize their results to $A_{\alpha}$-spectral radius as follows.

\begin{theorem}\label{thm3}
Let $T$ be an arbitrary tree on $n \geq 4$ vertices with the matching number $k$ such that $1\leq k\leq \lfloor\frac{n}{2}\rfloor$. If $0\leq\alpha<1$, then
$$\rho_{\alpha}(T)\leq \rho_{\alpha}(A(n,k))$$ with equality if and only if $T\cong A(n,k).$
\end{theorem}

\section{Proofs}

The first lemma is due to Nikiforov and Rojo [Lemma 6, \cite{NO}], it turns out to be very useful in the paper.
\begin{lemma}\label{le1}\rm\cite{NO}
Let $G$ be a connected graph with $\alpha\in [0,1)$. For $u,v\in V(G)$, suppose $N\subseteq N(v)\backslash (N(u)\cup \{u\})$. Let $G'=G-\{vw:w\in N\}+\{uw:w\in N\}$.
Let $X$ be a unit eigenvector of $A_{\alpha}(G)$ corresponding to $\rho_{\alpha}(G)$.
If $N\neq \emptyset$ and $x_u\geq x_v$, then $\rho_{\alpha}(G')>\rho_{\alpha}(G).$
\end{lemma}

\begin{lemma}\label{le2}\rm\cite{NO,X}
Let $G$ be a graph with a pendent path $P_{k}$. Suppose that the vertices in $P_{k}$ are labelled as $v_{1},v_{2},\cdots,v_{k}$ from outside to inside. Let $x$ be a unit eigenvector of $G$ corresponding to $\rho_{\alpha}(G)$. Denote by $x_{i}$ the entry of $x$ corresponding to vertices $v_{i}$. If $0\leq \alpha<1$ and $\rho_{\alpha}(G)>2$, then $x_{i}<x_{i+1}$ for $1\leq i\leq k-1$.
\end{lemma}

Before giving the proof of Theorem \ref{thm1}, we list all connected graphs with $A_{0}$-spectral radius at most 2 (see \cite{S}): the $A_{0}$-spectral radii of $H_{1}$, $H_{2}$, $H_{3}$, $H_{4}$ and $P_{n}$ are less than 2; the $A_{0}$-spectral radii of $H_{5}$, $H_{6}$, $H_{7}$, $H_{8}$ and $C_{k}$ are equal to 2. Let $G$ be a connected graph. If $G$ is not isomorphic to $P_{n}, C_{k}$ or $H_{i}$ for $i=1,2,3,4$, then \cite{VN1} Proposition 18 implies that $\rho_{\alpha}(G)>2$.

\newpage
\begin{center}
{\includegraphics{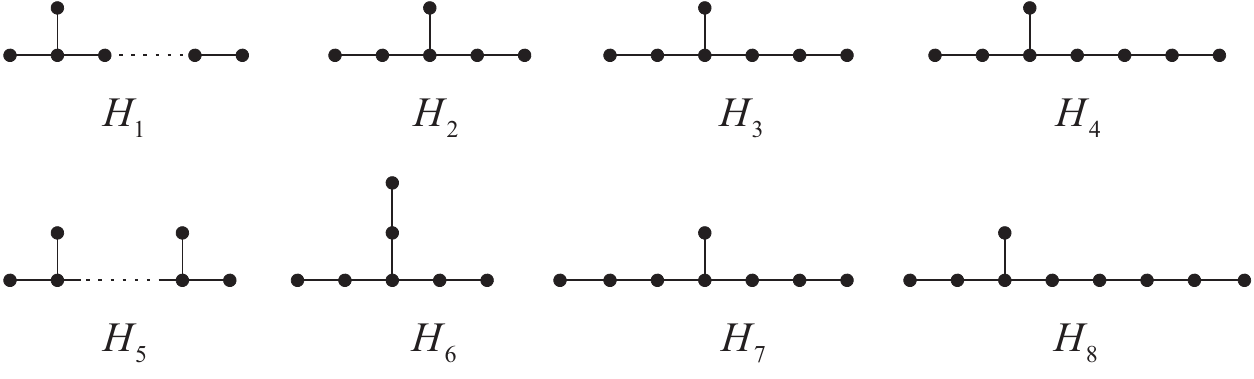}}
\vskip 0.2cm Fig. $1$. Graphs $H_{1}-H_{8}$.
\end{center}

\noindent{\bf{Proof of Theorem \ref{thm1}.}}
Let $G_{p,s,q}(u,v)$ be shown in Fig. 2 and $\rho_{\alpha}=\rho_{\alpha}(G_{p,s,q}(u,v))$. By contradiction, assume that $\rho_{\alpha}\geq \rho_{\alpha}(G_{p-1,s,q+1}(u,v))$. Suppose that $x$ is a Perron vector of $A_{\alpha}(G_{p,s,q}(u,v))$. For convenience, let $v=v_{q+1}$. The we have the following claim.
\begin{center}
{\includegraphics{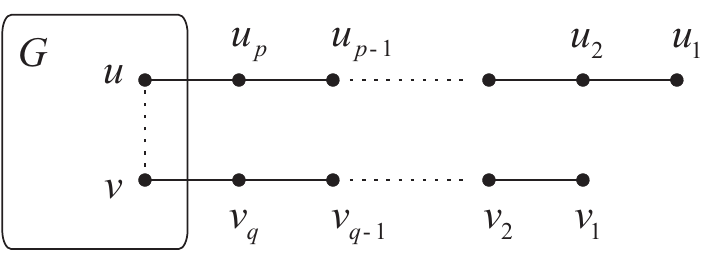}}
\vskip 0.2cm Fig. $2$. The graph $G_{p,s,q}(u,v)$.
\end{center}

\vspace{2mm}
\noindent{\bf Claim 1.} $x_{v_{i}}<x_{u_{i+1}}$ for $1\leq i\leq q+1$.
\vspace{2mm}

We will show this claim by induction. First, we shall show that $x_{v_{1}}<x_{u_{2}}$. If not, suppose that $x_{v_{1}}\geq x_{u_{2}}$ and let
$G_{1}=G_{p,s,q}(u,v)-u_{1}u_{2}+v_{1}u_{1}$ and it is clear that $G_{1}\cong G_{p-1,s,q+1}(u,v)$. According to Lemma \ref{le1}, we have $\rho_{\alpha}(G_{1})>\rho_{\alpha}$, a contradiction.
In the following, we assume that $x_{v_{i-1}}<x_{u_{i}}$ and we shall show that $x_{v_{i}}<x_{u_{i+1}}$ for $2\leq i\leq q+1$. If not, assume that $x_{v_{i}}\geq x_{u_{i+1}}$.
Let $G_{i}=G_{p,s,q}(u,v)-v_{i-1}v_{i}-u_{i}u_{i+1}+v_{i-1}u_{i+1}+v_{i}u_{i}$.
Clearly, $G_{i}\cong G_{p-1,s,q+1}(u,v)$. It follows that
\begin{eqnarray*}
\rho_{\alpha}(G_{i})-\rho_{\alpha}&\geq& x^{t}(A_{\alpha}(G_{i})-A_{\alpha}(G_{p,s,q}(u,v)))x\\
&=&2(1-\alpha)(x_{v_{i-1}}-x_{u_{i}})(x_{u_{i+1}}-x_{v_{i}})\\
&\geq& 0.
\end{eqnarray*}
Then $\rho_{\alpha}(G_{i})=\rho_{\alpha}$ and $x$ is also an eigenvector of $\rho_{\alpha}(G_{i})$. On the other hand,
$$0=\rho_{\alpha}(G_{i})x_{v_{i}}-\rho_{\alpha} x_{v_{i}}=[A_{\alpha}(G_{i})x]_{v_{i}}-[A_{\alpha}(G_{p,s,q}(u,v))x]_{v_{i}}=(1-\alpha)(x_{u_{i}}-x_{v_{i-1}})>0,$$ a contradiction. This completes the proof of the claim.

In the following, we will divide the left part of the proof into three cases.

\noindent{\bf{Case 1.}} $s=0$.

Let $G^\star=G_{p,0,q}(v,v)-\{wv|w\in N(v)\backslash \{v_{q},u_{p}\}\}+\{wu_{q+2}|w\in N(v)\backslash\{v_{q},u_{p}\}\}.$
It is clear that $G^\star\cong G_{p-1,0,q+1}(v,v)$. By Claim 1, we have $x_{v}=x_{v_{q+1}}<x_{u_{q+2}}$. Thus by Lemma \ref{le1}, we have $\rho_{\alpha}(G_{p-1,0,q+1}(v,v))>\rho_{\alpha}$, a contradiction.

\noindent{\bf{Case 2.}} $s=p-q-1$.

Let $G^{\star\star}=G_{p,s,q}(u,v)-\{wv|w\in {N(v)}\backslash \{v_{q},w_1\}\}+\{wu_{q+2}|w\in {N(v)}\backslash\{v_{q},w_1\}\}$. It is clear that $G^{\star\star}\cong G_{p-1,s,q+1}(u,v)$. By Claim 1, we have $x_{v}=x_{v_{q+1}}<x_{u_{q+2}}$. Thus by Lemma \ref{le1}, we have $\rho_{\alpha}(G_{p-1,s,q+1}(v,v))>\rho_{\alpha}$, a contradiction.

\noindent{\bf{Case 3.}}  $1\leq s\leq p-q-2$.

We shall first prove the following claim, for convenience, let $v=w_{0}$ and $u=w_{s}$.

\noindent{\bf Claim 2.} $x_{w_{i}}<x_{u_{q+2+i}}$ for $0\leq i\leq s$.

By induction. By Claim 1, we have $x_{w_{0}}=x_v<x_{u_{q+2}}$, so in the following we may assume that $x_{w_{i-1}}<x_{u_{q+1+i}}$.
Now let
\vspace{-3mm}
\begin{eqnarray*}
G_{q+1+i}&=&G_{p,s,q}(u,v)-\{wv|w\in {N(v)}\backslash \{v_{q},w_1\}\}-u_{q+2+i}u_{q+1+i}-w_{i-1}w_{i}\\
&&+\{wu_{q+2}|w\in {N(v)}\backslash\{v_{q},w_1\}\}+u_{q+1+i}w_{i}+u_{q+2+i}w_{i-1}.
\end{eqnarray*}
\begin{center}
{\includegraphics{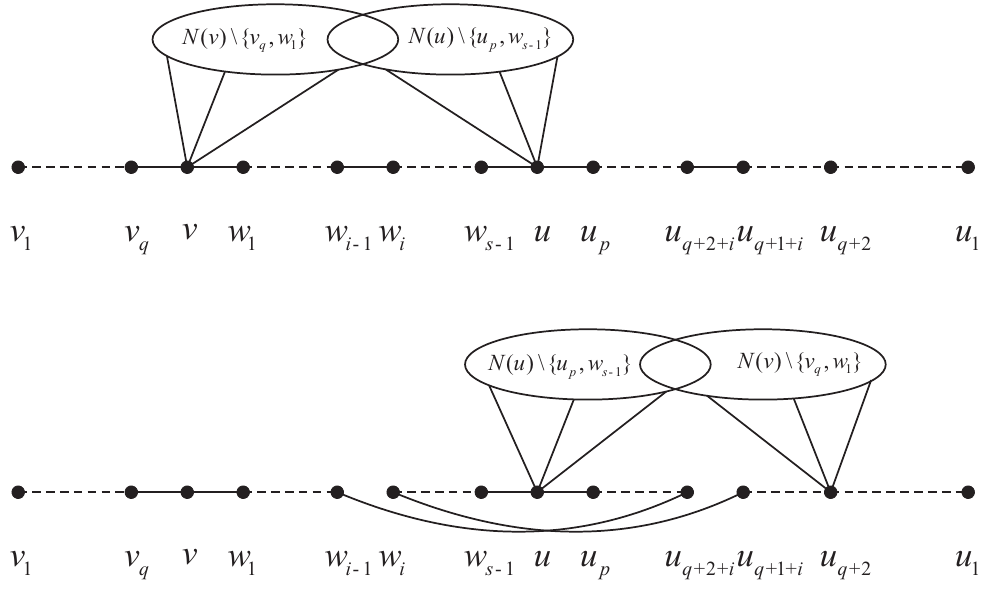}}
\vskip 0.2cm Fig. $3$. A transformation from $G_{p,s,q}(u,v)$ to $G_{p-1,s,q+1}(u,v)$.
\end{center}
Clearly, $G_{q+1+i}\cong G_{p-1,s,q+1}(u,v)$.
If $x_{w_i}\geq x_{u_{q+2+i}}$, then
\begin{eqnarray*}
&&\rho_{\alpha}(G_{q+1+i})-\rho_{\alpha}\geq x^{t}(A_{\alpha}(G_{2})-A_{\alpha}(G_{p,s,q}(u,v)))x\\
&=&2(1-\alpha)(x_{w_{i}}-x_{u_{q+2+i}})(x_{u_{q+1+i}}-x_{w_{i-1}})+\alpha(d(v)-2)(x_{u_{q+2}}^2-x_{w_{0}}^2)\\
&&+2(1-\alpha)(x_{u_{q+2}}-x_{w_{0}})\sum_{w\in {N(v)}\backslash \{v_{q},w_1\}}x_w\\
&>& 0,
\end{eqnarray*}
a contradiction. Therefore, we have $x_{w_i}<x_{u_{q+2+i}}$, which completes the proof of this claim.

According to Claim 2, one can see that $x_{u}=x_{w_{s}}<x_{u_{q+2+s}}$ since $p\geq q+s+2$.
On the other hand, note that $G_{p,s,q}(u,v)$ contains $H_5$ as a proper induced subgraph.
Then $\rho_{\alpha}> \rho_{\alpha}(H_5)\ge2.$
By Lemma \ref{le2}, we have $x_{u}>x_{u_{p}}\geq x_{u_{q+2+s}}$, a contradiction.
This completes the proof.  {\hfill$\Box$}

\vspace{3mm}
\noindent\textbf{Proof of Theorem \ref{thm2}.}
Suppose that $G$ attains the maximum $A_{\alpha}$-spectral radius among all connected graphs with order $n$ and $k$ cut vertices.
Since $A_{\alpha}(G)$ is irreducible, we have $\rho_{\alpha}(G)<\rho_{\alpha}(G+e)$ if $e\notin E(G)$.
Consequently, we can assume that each cut vertex of $G$ connects exactly two blocks and that all of these blocks are cliques.
Assume that all the blocks of $G$ are $K_{a_1},K_{a_2},\ldots, K_{a_{k+1}}$. Then we have the following claim.

\vspace{3mm}
\noindent \textbf{Claim.} There is no cut vertex belongs to $K_{a_{i}}$ and $K_{a_{j}}$ with $a_{i}\geq 3, a_{j}\geq 3$.

If not, suppose that $v$ is the cut vertex which belongs to $K_{a_{i}}$ and $K_{a_{j}}$. Without loss of generality, we may assume that $u$ is a vertex in $V(K_{a_{i}})$ such that
$x_{u}=\max\{x_{w}|w\in V(K_{a_{i}})\cup V(K_{a_{j}})\backslash \{v\}\}$. Let $u'$ be a vertex in $ V(K_{a_{j}})\backslash \{v\}$.
Let $G'=G-\{wu'|w\in V(K_{a_{j}})\backslash \{v,u'\}\}+\{wu|w\in V(K_{a_{j}})\backslash \{v,u'\}\}$. Note that $G$ and $G'$ have the same number of cut vertices and $x_{u}\geq x_{u'}$. Then by Lemma \ref{le1}, we have $\rho_{\alpha}(G')>\rho_{\alpha}(G)$, which contradicts the maximality of $\rho_{\alpha}(G)$.

Ordering the cardinalities of these blocks as $a_1 \geq a_2 \geq\cdots\geq a_{k+1} \geq 2.$
In the following, we shall show that $G$ has exactly one clique with order greater than 2. If not, we assume that $a_{1}\geq a_{2}\geq 3$.
Hence, $a_{1}=n+k-\sum_{i=2}^{k+1}a_{i}\leq n-k-1$. Then by the \textbf{Claim}, we know that $\Delta\leq a_{1}\leq n-k-1$.
Moreover, since $G$ is clearly irregular, it follows that $\rho_{\alpha}(G)<\Delta\leq n-k-1$. On the other hand, $\rho_{\alpha}(G_{n,k})>\rho_{\alpha}(K_{n-k})=n-k-1$. Therefore, $\rho_{\alpha}(G)<\rho_{\alpha}(G_{n,k})$, contradicting the maximality of $\rho_{\alpha}(G)$. Thus, it follows that $a_{1}\geq 2=a_{2}=\cdots=a_{k+1}$, that is, $G$ is a graph formed by attaching some pendent paths to a clique $K_{n-k}$. Now, by Theorem \ref{thm1}, we obtain $G\cong G_{n,k}$.  {\hfill$\Box$}

In order to give the proof of Theorem \ref{thm3}, the following transformation is needed, which is a direct result from Theorem \ref{thm1}.

\noindent{\bf{Transformation A}} Let $T_1$ be a tree with $n-p$ vertices and $v$ be a vertex of $V(T_1)$.
Let $T$ be the tree by attaching a pendent path $v$ of $T_1$ with length $p$ and $T'$ be the tree by attaching two pendent paths at $v$ of $T_1$, with lengths $2$ and $p-2$, respectively.
Then $m(T)=m(T')$ and $\rho_{\alpha}(T')>\rho_{\alpha}(T).$

\vspace{3mm}
\noindent\textbf{Proof of Theorem \ref{thm3}.}
Suppose that $T$ attains the maximum $A_{\alpha}$ spectral radius among all trees with matching number $m.$
Then by {\bf{Transformation A}}, the length of every external path of $T$ is at most 3. Let $x$ be the Perron vector
of $\rho_{\alpha}(T)$ and $u$ be a vertex with $x_u=\max\{x_w|w\in V(T)\}$. In the following, we only need to show that there
are no vertices of $V(T)$ with degree greater than 2 except for $u$. If not, suppose that $d(v)\geq 3$ and $d(u,v)$ is maximal.
Let $N(v)=\{v_1\}\cup V_1\cup V_2$ where $v_1$ is the vertex in the $(u,v)$-path, $V_1=\{w|d(w)=1\}$ and $V_2=\{w|d(w)=2, w\neq v_1\}$.
If $|V_1|=\varnothing$, then let $T'=T-\{vw|w\in V_2\}+\{uw|w\in V_2\}$. Note that $m(T)=m(T')$ and $x_u\geq x_v$, then by Lemma \ref{le1},
we have $\rho_{\alpha}(T')>\rho_{\alpha}(T)$, a contradiction. In the following, we assume that $|V_1|\neq\varnothing$. Then either $vv_1$ or $vw$ ($w\in V_1$) is in
a maximal matching. If attaching a pendent vertex to $u$ will not increase the matching number, then let $T'=T-\{vw|w\in V_1\cup V_2\backslash\{s\}, s\in V_1\}+\{uw|w\in V_1\cup V_2\backslash\{s\}, s\in V_1\}$. Note that $m(T)=m(T')$ and $x_u\geq x_v$, then again by Lemma \ref{le1},
we have $\rho_{\alpha}(T')>\rho_{\alpha}(T)$, a contradiction.
If attaching a pendent vertex to $u$ will increase the matching number, then let $T'=T-\{vw|w\in V_1\cup V_2\}+\{uw|w\in V_1\cup V_2\}-vv_1+uv$.
It is clear that $m(T)=m(T')$ and
\begin{eqnarray*}
\rho_{\alpha}(T')-\rho_{\alpha}(T)&\geq& \alpha d(v)x_u-\alpha (d(v)-1)x_v-\alpha x_{v_1}\\
&+&(1-\alpha)[\sum_{w\in V_1\cup V_2}x_w(x_u-x_v)+x_v(x_u-x_{v_1})]\\
&\geq&0.
\end{eqnarray*}
If $\rho_{\alpha}(T')=\rho_{\alpha}(T)$, then $x$ is also the Perron vector of $A_{\alpha}(T')$.
On the other hand, $\rho_{\alpha}(T')x_u=\rho_{\alpha}(T)x_u+\sum_{w\in V_1\cup V_2}x_w>\rho_{\alpha}(T)x_u,$ a contradiction.
Thus we have $\rho_{\alpha}(T')>\rho_{\alpha}(T)$, which contradicts the fact that $\rho_{\alpha}(T)$ is maximum. This completes the proof. {\hfill$\Box$}

\end{document}